\documentclass[a4paper,secthm,amsthm]{elsart}
\usepackage{amsmath,amssymb,amsfonts,euscript} 

\newlength{\textequation} 
\DeclareMathOperator{\st}{|} 
\DeclareMathOperator{\otensprod}{\otimes} 
\DeclareMathOperator{\ctensprod}{\overline{\otimes}}
\DeclareMathOperator{\spec}{\sigma} 

\newcommand{\C}{\mathbb{C}}
\newcommand{\R}{\mathbb{R}}
\newcommand{\N}{\mathbb{N}}

\providecommand{\abs}[1]{\lvert #1 \rvert}
\providecommand{\scal}[2]{\langle #1,#2 \rangle}
\providecommand{\norm}[1]{\lVert #1 \rVert}
\providecommand{\set}[1]{\left\{ #1 \right\}}

\newcommand{\cd}[1]{\mathop{\bf #1}\nolimits}

\newcommand{\res}[1]{_{\rceil #1}}
\newcommand{\adh}[1]{\, \overline{#1}}

\renewcommand{\geq}{\geqslant}
\renewcommand{\leq}{\leqslant}

\def\eps{\varepsilon}

\def\vv{\underline{v}}

\def\zz{\underline{\zeta}}

\def\kk{\underline{k}}

\def\kapt{\tilde{\kappa}}

\def\xt{\tilde{x}}

\def\et{\tilde{e}}

\def\AAt{\tilde{\A}}
\def\At{\tilde{A}}
\def\Xt{\tilde{X}}
\def\Ht{\tilde{H}}
\def\Mt{\tilde{M}}
\def\nt{\tilde{n}}

\def\comp{_{\text{comp}}}
\def\loc{_{\text{loc}}}

\def\Tmin{T_{\text{min}}} 

\def\Deltat{\tilde{\Delta}}
\def\A{\mathcal{A}}
\def\B{\mathcal{B}}
\def\CC{\mathcal{C}}
\def\sA{\sqrt{\!  A}}
\def\sAp{\sqrt{\!  A'}} 

\def\Ap{A} 
\def\Am{A'} 
\def\As{\A'{}^{s}}
\def\AT{A_{\theta}}

\def\ATb{A_{-\theta}}

\providecommand{\noz}[1]{\norm{#1}_{0}}
\providecommand{\noo}[1]{\norm{#1}_{1}}
\providecommand{\nod}[1]{\norm{#1}_{2}}
\providecommand{\nom}[1]{\norm{#1}_{-1}}
\providecommand{\nop}[1]{\norm{#1}_{p}}
\providecommand{\nopm}[1]{\norm{#1}_{p-1}}
\providecommand{\nomp}[1]{\norm{#1}_{-p}}
\providecommand{\scd}[2]{\scal{#1}{#2}}
\providecommand{\scz}[2]{\scal{#1}{#2}_{0}}

\def\dx{\dot{x}}
\def\dxi{\dot{\xi}}
\def\dz{\dot{z}}
\def\ddz{\ddot{z}}
\def\dza{\dot{\zeta}}
\def\ddza{\ddot{\zeta}}
\def\dzz{\dot{\zz}}
\def\ddzz{\ddot{\zz}}

\def\dv{\dot{v}}
\def\dvv{\dot{\vv}}
\def\df{\dot{f}}

\def\dphi{\dot{\phi}}
\def\dchi{\dot{\chi}}
\def\dchiT{\dot{\chi}_{T}}
\def\dS{\dot{S}}

\def\L{\mathcal{L}} 
\def\D{D} 

\begin{document}


\begin{frontmatter}
\title{Controllability cost of conservative systems:\\
resolvent condition and transmutation.}
\author{Luc Miller}
\address{
{\'E}quipe Modal'X, EA 3454, \\
Universit{\'e} Paris X, B{\^a}t.~G,
200 Av.~de la R{\'e}publique,
92001 Nanterre, France.}
\address{Centre de Math{\'e}matiques Laurent Schwartz, UMR CNRS 7640, \\ 
{\'E}cole Polytechnique, 91128
Palaiseau, France.}
\ead{miller@math.polytechnique.fr}

\thanks{This work was partially supported by the ACI grant
  ``{\'E}quation des ondes : oscillations, dispersion et contr{\^o}le''.
First on {\tt arXiv:math.OC/0402058}, 4 February 2004.
Accepted by the Journal of Functional Analysis 10 February 2004.}

\begin{abstract} 
This article concerns the exact controllability of unitary groups 
on Hilbert spaces with unbounded control operator. 
It provides a necessary and sufficient condition not involving time
which blends a resolvent estimate and an observability inequality. 
By the transmutation of controls in some time $L$
for the corresponding second order conservative system,
it is proved that the cost of controls in time $T$ for the  unitary group 
grows at most like $\exp(\alpha L^{2}/T)$
as $T$ tends to $0$. 
In the application to the cost of fast controls 
for the Schr{\"o}dinger equation,
$L$ is the length of the longest ray of geometric optics
which does not intersect the control region. 
This article also provides observability resolvent estimates implying
fast smoothing effect controllability at low cost,
and underscores that the controllability cost of a system 
is not changed by taking its tensor product with a conservative system. 
\end{abstract}






\end{frontmatter}
\thispagestyle{headings} 


\section{Introduction}
\label{sec:intro}

Let $H_{0}$ and $Y$ be Hilbert spaces 
with respective norms $\noz{\cdot}$ and $\norm{\cdot}$.
Let $A$ be a self-adjoint, positive and boundedly invertible 
unbounded operator on $H_{0}$ with domain $\D(A)$.
We introduce the Sobolev scale of spaces based on $A$.
For any positive integer $p$, 
let $H_{p}$ denote the Hilbert space $\D(A^{p/2})$ 
with the norm $\nop{x}=\noz{A^{p/2} x}$
(which is equivalent to the graph norm $\noz{x}+\noz{A^{p/2} x}$).
We identify $H_{0}$ and $Y$ with their duals.
Let $H_{-p}$ denote the dual of $H_{p}$.
Since $H_{p}$ is densely continuously embedded in $H_{0}$,
the pivot space $H_{0}$ is densely continuously embedded in $H_{-p}$,
and $H_{-p}$ is the completion of $H_{0}$ with respect to the norm 
$\nomp{x}=\noz{A^{-p/2} x}$.
We still denote by $\Ap$
the restriction of $A$ to $H_{p}$ with domain $H_{p+2}$. 
It is self-adjoint with respect to the $H_{p}$ scalar product.
We denote by $A'$ its dual 
with respect to the duality between $H_{p}$ and $H_{-p}$,
which is an extension  of $A$ to $H_{-p}$ with domain $H_{2-p}$.

Let $C\in \L(H_{2};Y)$ 
and let $B\in\L(Y,H_{-2})$ denote the dual of $C$
(where $\L(X,Y)$ denotes 
the Banach space of continuous operators from $X$ to $Y$).

The ``~generator~'' $\Am$ and the control operator $B$ 
define the first and second order differential equations:
\begin{gather}\label{eqfsystcon}
\dphi(t)-i\Am\phi (t)=Bu(t), \quad 
\phi(0)=\phi_{0}\in H_{-1}, \quad
u\in L^{2}\loc(\R;Y) , \\
\label{eqsystcon}
\ddza(t)+\Am\zeta (t)=Bv(t),  
\zeta(0)=\zeta_{0}\in H_{0}, \dza(0)=\zeta_{1}\in H_{-1}, v\in L^{2}\loc(\R;Y) 
, \end{gather}
where
each dot denotes a derivative with respect to the time variable $t$,
$u$ and $v$ are the input functions.

Equations  (\ref{eqfsystcon}) and (\ref{eqsystcon})   with $u=0=v$ 
describe reversible conservative systems.
For example, if $A$ is the positive Laplacian and 
$B$ is a boundary control operator,
then (\ref{eqsystcon}) is a boundary controlled scalar wave equation,
(\ref{eqfsystcon}) is 
a boundary controlled Schr{\"o}dinger equation 
(sections~\ref{sec:1d} and~\ref{sec:appli} elaborate on this example).

We assume that $B$ is 
{\em an admissible control operator for (\ref{eqsystcon})}, i.e.:
\begin{gather}
\label{eqadm}
\forall T>0, \forall v\in L^{2}(0,T;Y), \quad
\int_{0}^{T}e^{it\sAp}Bv(t)dt \in H_{-1}
, \end{gather} 
so that the solution $\zeta\in C^{0}(\R;H_{0})\cap C^{1}(\R;H_{-1})$ 
of (\ref{eqsystcon})
is defined by the following integral formula
where $S(t)=(\sAp)^{-1}\sin(t\sAp)$ and $\dS(t)=\cos(t\sAp)$:
\begin{gather*}
\zeta(t)=\dS(t)\zeta_{0}+S(t)\zeta_{1}+\int_{0}^{t}S(t-s)Bv(s)ds
. \end{gather*}
In section~\ref{sec:setting} the control system (\ref{eqsystcon})
and its dual observation system are reduced to 
the standard first-order setting for the theory of observation and control. 
We also assume that $B$ is 
{\em an admissible control operator for (\ref{eqfsystcon})}, 
i.e.: $\forall T>0$, $\exists K_{1,T}>0$,
\begin{gather}
\label{eqfadm}
\forall u\in L^{2}(0,T;Y), \quad
\nom{\int_{0}^{T}e^{it\Am}Bu(t)dt}^{2} 
\leq K_{1,T} \int_{0}^{T}\norm{u(t)}^{2}dt
, \end{gather} 
so that the solution
$\phi\in C^{0}(\R;H_{-1})$ of (\ref{eqfsystcon})
is defined by the integral formula:
\begin{gather}\label{eqfintcon}
\phi(t)=e^{it\Am}\phi_{0}+\int_{0}^{t}e^{i(t-s)\Am}Bu(s)ds
. \end{gather}

\begin{defn}
\label{def:con}
The system (\ref{eqfsystcon}) 
is {\em exactly controllable in time $T$} if 
for all $\phi_{0}$ in $H_{-1}$, 
there is a $u$ in $L^{2}(\R;Y)$ such that 
$u(t)=0$ for $t\notin [0,T]$ and $\phi(T)=0$.
The {\em controllability cost} for (\ref{eqfsystcon}) in time $T$ is
the smallest positive constant $\kappa_{1,T}$ 
in the following inequality for all such $\phi_{0}$ and $u$:
\begin{gather}
\label{eqkappa1} 
\int_{0}^{T}\norm{u(t)}^{2}dt\leq \kappa_{1,T}\nom{\phi_{0}}^{2}
. \end{gather} 
The system (\ref{eqsystcon}) 
is {\em exactly controllable in time $T$} if 
for all $\zeta_{0}$ in  $H_{0}$ 
and $\zeta_{1}$ in  $H_{-1}$, 
there is a $v$ in $L^{2}(\R;Y)$ such that 
$v(t)=0$ for $t\notin [0,T]$ and 
$\zeta(T)=\dza(T)=0$.
The {\em controllability cost} for (\ref{eqsystcon}) in time $T$ is
the smallest positive constant $\kappa_{2,T}$ 
in the following inequality for all such 
$\zeta_{0}$, $\zeta_{1}$ and $v$:
\begin{gather}
\label{eqkappa2} 
\int_{0}^{T}\norm{v(t)}^{2}dt\leq \kappa_{2,T}
\left( \noz{\zeta_{0}}^{2}+\nom{\zeta_{1}}^{2} \right)
. \end{gather} 
\end{defn}

\begin{rem}
Strictly speaking, 
the properties above define the null-controllability of the systems
but, for such reversible systems, they are equivalent to exact controllability.
We refer to section~\ref{sec:setting} for the dual notions 
of observability.
\end{rem}

The main results of the paper, stated in section~\ref{sec:results}, 
are consequences 
of the controllability of the wave-like system (\ref{eqsystcon}) 
on the controllability of the Schr{\"o}dinger-like system (\ref{eqfsystcon}).
In particular, upper bounds on the 
controllability cost $\kappa_{1,T}$
of (\ref{eqfsystcon}) as $T$ tends to zero are given.
Applications to 
the boundary controllability of the Schr{\"o}dinger equation,
based on the geodesic condition of Bardos-Lebeau-Rauch (cf.~\cite{BLR92})
for the controllability of the wave equation,
are presented in section~\ref{sec:appli}:
a new proof of the result of Lebeau in~\cite{Leb92}, 
an extension of this result to product manifolds,
and an upper bound on the cost of fast controls in the same context.

The main tool presented in this paper is 
the {\em control transmutation method} 
which can be seen as an adaptation to the theory of control 
of the kernel estimates method of Cheeger-Gromov-Taylor in \cite{CGT82}.
It consists in explicitly constructing controls $v$
in any time $T$ for the Schr{\"o}dinger-like system (\ref{eqfsystcon})
in terms of controls $u$ in time $L$
for the corresponding wave-like system (\ref{eqsystcon}),
i.e. $u(t,x)=\int_{\mathbb{R}} k(t,s)v(s)\, ds $,
where the compactly supported kernel $k$ which depends on $T$ and $L$ is 
some {\em fundamental controlled solution} on the segment $[-L,L]$
controlled at both ends.
In section~\ref{sec:1d}, we recall an earlier 
estimate on the optimal fast control cost rate 
for a one dimensional system.
We use it to construct the fundamental controlled solution $k$
in section~\ref{sec:fcs}
and perform the transmutation in section~\ref{sec:transmut}.

This paper also contains results of independent interest 
on the controllability of systems defined by unitary groups.
In section~\ref{sec:setting}, 
we recall admissibility, observability and controllability
notions for such systems, their duality, 
and reduce the second-order system  (\ref{eqsystcon}) 
to this first-order setting. 
In section~\ref{sec:res}, we state a necessary and sufficient condition 
on the resolvent of the generator and the observation operator  
for exact observability. 
In section~\ref{sec:smooth}, we state a sufficient condition 
on the resolvent of the generator and the observation operator 
for the existence of controls steering any state to a smooth state
in any positive time $T$ 
at a cost bounded from above by a negative power of $T$.
In section~\ref{sec:tensor}, we prove that 
the controllability cost of a system 
is not changed by taking its tensor product with a system
defined by a unitary group.

In the companion paper~\cite{LMctm},
we apply the control transmutation method 
to the simpler case of the first-order equation 
$\dphi(t)+e^{i\theta}\Am\phi (t)=Bu(t)$ 
with $\abs{\theta}<\pi/2$
(in particular, the correponding semigroup is holomorphic),
but in the more general setting 
where $A'$ generates a cosine operator function in a Banach space.
The relationship between the controllability of 
this first order equation and the second order equation (\ref{eqsystcon})
has been investigated earlier in various settings
(cf. references in~\cite{LMctm}) but only with $\theta=0$.
Equation (\ref{eqfsystcon}) corresponds to the case $\abs{\theta}=\pi/2$.


\section{Boundary control of the Schr{\"o}dinger equation on a segment}
\label{sec:1d}

Our estimate of the cost of fast controls for (\ref{eqfsystcon}) builds,
through the control transmutation method, 
on the same estimate 
for a one dimensional control system of type (\ref{eqfsystcon}),
i.e. the Schr{\"o}dinger equation on a segment $[0,L]$
with Dirichlet ($N=0$) or Neumann ($N=1$) condition at the left end
controlled at the right end through a Dirichlet condition: 
\begin{gather}
\label{eq1dp}
\partial_{t}\phi + i\partial_{s}^{2} \phi=0
\ \text{on}\ \left]0,T\right[\times \left]0,L\right[, \ 
\partial_{s}^{N}\phi\res{s=0} = 0, \ 
\phi\res{s=L}=
u,\
\phi\res{t=0} = \phi_{0}
. \end{gather}
With the notations of section~\ref{sec:intro}, 
$A=-\partial_{s}^{2}$ on $H_{0}=L^{2}(0,L)$
with $\D(A)=\{f\in H^{2}(0,L) \,|\, \partial_{s}^{N}f(0)=f(L)=0\}$,
$C$ with values in $Y=\C$ is defined by $Cf=\partial_{s}f(L)$, 
and $H_{1}=H^{1}_{N}$ is one of the following 
Sobolev spaces on the segment $[0,L]$: 
\begin{gather*}
H^{1}_{1}(0,L)=\{f\in H^{1}(0,L) \,|\, f(L)=0\}
 \mbox{ and }  
H^{1}_{0}(0,L)=\{f\in H^{1}_{1}(0,L) \,|\, f(0) = 0\}
\ .
\end{gather*}
\begin{defn}
\label{defin:alpha}
The rate $\alpha_{*}$ is the smallest positive constant such that 
for all $\alpha > \alpha_{*}$
there exists $\gamma>0$ such that, 
for all $N\in\{0,1\}$,
$L>0$, $T\in \left]0,\inf(\pi,L)^{2}\right]$
the controllability cost $\kappa_{L,T}$ of the system (\ref{eq1dp})
satisfies: $\kappa_{L,T}\leq \gamma \exp(\alpha L^{2}/T)$.
\end{defn}
It is well-known that the controllability of this system 
reduces by spectral analysis to classical results 
on nonharmonic Fourier series.
The study of upper bounds of the controllability cost for short times 
was initiated by Seidman (cf. references in \cite{LMschrocost}). 
We recall a theorem of~\cite{LMschrocost} 
which improves his estimate of the optimal rate $\alpha_{*}$
(computing $\alpha_{*}$ is an interesting open problem and its solution 
does not have to rely on the analysis of series of complex exponentials).
\begin{thm}
\label{th:1d}
The optimal fast control cost rate for the one dimensional system (\ref{eq1dp})
in definition~\ref{defin:alpha} 
satisfies: $1/2\leq \alpha_{*}\leq 8\left(36/37\right)^{2}< 8$.
\end{thm}


\section{Main results}
\label{sec:results}

The main result of this paper, proved in section~\ref{sec:transmut}, 
is a generalization of theorem~\ref{th:1d} to 
the first-order system (\ref{eqfsystcon}) 
under some condition on the second-order system (\ref{eqsystcon}):
\begin{thm}
\label{th}
If the system (\ref{eqsystcon}) is exactly controllable 
for times greater than $L_{*}$,
then the system (\ref{eqfsystcon}) is exactly controllable 
in any time $T$.
Moreover, the controllability cost $\kappa_{1,T}$ of (\ref{eqfsystcon}) 
satisfies the following upper bound 
(with $\alpha_{*}$ as in theorem~\ref{th:1d}):
\begin{gather}
\label{equb}
\limsup_{ T\to 0}  
T \ln \kappa_{1,T} \leq \alpha_{*} L_{*}^{2}
\ .\end{gather}
\end{thm}

\begin{rem}
\label{rem:th}
The upper bound (\ref{equb}) means that 
the smallest norm of an input function $u$
steering the system (\ref{eqfsystcon}) from an initial state $\phi_{0}$
to zero grows at most like $\gamma\norm{\phi_{0}}\exp(\alpha L_{*}^{2}/(2T))$
as the control time $T$ tends to zero 
(with any $\alpha>\alpha_{*}$ and some $\gamma>0$).
The falsity of the converse of theorem~\ref{th} is well-known, 
e.g. in the setting of section~\ref{sec:appli}. 
\end{rem}
\begin{rem}
As observed in~\cite{Car93} 
(\ref{equb}) yields a
logarithmic modulus of continuity for the minimal time function 
$\Tmin:H_{-1}\to [0,+\infty)$ of (\ref{eqfsystcon}); i.e.
$\Tmin(\phi_{0})$, defined as the infimum of the times $T>0$ for which 
there is a $u$ in $L^{2}(\R;Y)$ such that 
$\int_{0}^{T}\norm{u(t)}^{2}dt\leq 1$,
$u(t)=0$ for $t\notin [0,T]$ and $\phi(T)=0$,
satisfies: for all $\alpha>\alpha_{*}$, 
there is a $c>0$ such that, 
for all $\phi_{0}$ and $\phi_{0}'$ in $H_{-1}$
with $\nom{\phi_{0}-\phi_{0}'}$ small enough,
$\abs{\Tmin(\phi_{0}) - \Tmin(\phi_{0}')}\leq 
\alpha L_{*}^{2}/\ln(c/\nom{\phi_{0}-\phi_{0}'})$.
\end{rem}

\medskip

Replacing the notion of exact controllability 
by the controllability to a subspace with finite spectrum,
which is enough to steer any initial state to a smooth final state,
we obtain a much better upper bound for the cost of fast controls.
The spectral projection on $[\lambda_{1}, \lambda_{2}]$
is denoted by $\cd{1}_{\lambda_{1}\leq A\leq \lambda_{2}}$.

\begin{thm}
\label{th:smoothing}
If the system (\ref{eqsystcon}) is exactly controllable,
then $\exists \kappa>0$, $\exists d>0$, 
$\forall T\in ]0,1]$, $\forall \phi_{0} \in H_{-1}$, 
$\exists u \in L^{2}(\R;Y)$ such that 
the solution $\phi\in C^{0}([0,\infty);H_{-1})$ of (\ref{eqfsystcon})
satisfies $\cd{1}_{\abs{A}\geq d/T^2}\phi(T)=0$
and 
$\int_{0}^{T}\norm{u(t)}^{2}dt\leq 
\frac{\kappa}{T}\nom{\cd{1}_{\abs{A}\geq d/T^2} \phi_{0}}^{2}$.
In particular, for all $p\in \N$: 
$\phi(T)\in H_{p-1}$ and 
$\nopm{\phi(T)}
\leq (1+\sqrt{\kappa K_{1,T}/T})(d/T^2)^{p/2} \nom{\phi_{0}}$.
\end{thm}

\medskip

Theorem~\ref{th} 
still holds when the system (\ref{eqfsystcon}) 
is replaced by its tensor product with a conservative system.
If we consider $A$ as a self-adjoint operator on $H_{1}$ 
and if $\At$ is an other self-adjoint operator on an other Hilbert space $\Ht$,
then the operator $A\otensprod I+I\otensprod\At$ defined on 
the algebraic tensor product $D(A)\otensprod D(\At)$ is closable
and its closure, denoted $A+\At$, 
is a self-adjoint operator on 
the closure of the algebraic tensor products $H_{1}\otensprod \Ht$,
denoted $H_{1}\ctensprod \Ht$ (cf. theorem~VIII.33 in~\cite{RS}).
The self-adjoint operator $A'+\At$ is defined similarly.
Thanks to lemma~\ref{lem:abstprod} proved in section~\ref{sec:tensor}
(and the duality between observability and controllability),
theorem~\ref{th} implies:

\begin{thm}
\label{th:prod}
Let $\At$ be a self-adjoint operator on an other Hilbert space $\Ht$.
If the system (\ref{eqsystcon}) is exactly controllable 
in times greater than $L_{*}$,
then for all positive time $T$ there is a positive constant $\kapt_{T}$ 
satisfying (with $\alpha_{*}$ as in theorem~\ref{th:1d}):
\begin{gather*}
\forall F\in H_{1}\ctensprod \Ht, \ 
\int_{0}^{T}\norm{(C\otensprod I)e^{it(A+\At)}F}^{2} dt
\leq \kapt_{T} \norm{F}^{2}
\ \text{and} \ \limsup_{ T\to 0}  T \ln \kapt_{T} \leq \alpha_{*} L_{*}^{2}
\, .\end{gather*}
This is equivalent to the exact controllability in time $T$ at cost $\kapt_{T}$
of the equation $\dot{\Phi}(t)-i(A'+\At)\Phi(t)=(B\otimes I)u(t)$
with $\Phi(0)=\Phi_{0}\in H_{-1}\ctensprod \Ht$ and $u\in L^{2}\loc(\R;Y\ctensprod \Ht)$.
\end{thm}


\section{Preliminaries on conservative control systems}
\label{sec:setting}

In this section, we review
the general setting for conservative control systems:
admissibility, observability and controllability
notions and their duality (cf.~\cite{DR77} and~\cite{Wei89}).
We recall the characterization of solutions in the weak sense.
We prove that smoother data can be controlled with smoother input functions.
We reduce the second-order system  (\ref{eqsystcon})
to this first-order setting.

Let $X$ and $Y$ be Hilbert spaces.
Let $\A:\D(\A)\to X$ be a self-adjoint operator.
Equivalently, $i\A$ generates a strongly continuous group $(e^{it\A})_{t\in\R}$ 
of unitary operators on $X$. 
Let $X_{1}$ denote $\D(\A)$ with the norm $\norm{x}_{1}=\norm{(\A-\beta)x}$
for some $\beta\notin \spec(\A)$
($\spec(\A)$ denotes the spectrum of $\A$, 
this norm is equivalent to the graph norm 
and $X_{1}$ is densely and continuously embedded in $X$)
and let $X_{-1}$ be the completion of $X$ with respect to the norm 
$\norm{\xi}_{-1}=\norm{(\A-\beta)^{-1}\xi}$.
Let $X'$ denote the dual of $X$ 
with respect to the pairing $\scd{\cdot}{\cdot}$
(linear in the first variable and conjugate-linear in the second variable). 
The dual of $\A$ is a self-adjoint operator $\A'$ on $X'$.
The dual of $X_{1}$ is the space $X'_{-1}$ 
which is the completion of $X'$ with respect to the norm 
$\norm{\xi}_{-1}=\norm{(\A'-\bar{\beta})^{-1}\xi}$
and the dual of $X_{-1}$ is the space $X'_{1}$
which is $\D(\A')$ 
with the norm $\norm{x}_{1}=\norm{(\A'-\bar{\beta})x}$.

Let $\CC\in \L(X_{1},Y)$ and let $\B\in\L(Y',X'_{-1})$ denote its dual.
Note that the same theory applies to any $\A$-bounded operator $\CC$ 
with a domain invariant by $(e^{it\A})_{t\geq 0}$ 
since it can be represented by an operator in $\L(X_{1},Y)$ (cf.~\cite{Wei89}).

We consider the dual observation and control systems
with output function $y$ and input function $u$:
\begin{gather}
\dx(t)-i\A x(t)=0, \quad x(0)=x_{0}\in X, 
\quad y(t)=\CC x(t), \label{eqsystx}\\
\dxi(t)-i\A' \xi(t)=\B u(t), \quad \xi(0)=\xi_{0}\in X', 
\quad u\in L^{2}\loc(\R;Y') \label{eqsystxi}
. \end{gather}
  
We make the following equivalent admissibility assumptions 
on the observation operator $\CC$ and the control operator $\B$
(cf.~\cite{Wei89}): $\forall T>0$, $\exists K_{T}>0$,
\begin{gather}
\label{eqadmobs}
\forall x_{0}\in \D(\A),
\quad
\int_{0}^{T}\norm{\CC e^{it\A}x_{0}}^{2}dt \leq K_{T}\norm{x_{0}}^{2} , \\
\label{eqadmcon}
\forall u\in L^{2}(\R;Y'), \quad
\norm{\int_{0}^{T}e^{it\A'}\B u(t)dt}^{2}
\leq K_{T} \int_{0}^{T}\norm{u(t)}^{2}dt .
\end{gather} 
With this assumption, the output map 
$x_{0} \mapsto y$ 
from $\D(\A)$ to $L^{2}\loc(\R;Y)$ 
has a continuous extension to $X$. 
The equations (\ref{eqsystx}) and (\ref{eqsystxi}) 
have unique solutions $x\in C(\R,X)$ and $\xi\in C(\R,X')$ 
defined by:
\begin{gather}
\label{eqmildsol}
x(t)=e^{it\A}x_{0},\quad 
\xi(t)=e^{it\A'}\xi(0)+\int_{0}^{t}e^{i(t-s)\A}\B u(s) ds .
\end{gather}
These so-called mild solutions are also 
the unique solutions in the weak sense (cf. \cite{Bal77}):
$x(0)=x_{0}$, $\xi(0)=\xi_{0}$,  
\begin{gather}
\label{eqweaksolx}
\forall \varphi\in \D(\A'), \ t\mapsto \scd{x(t)}{\varphi}\in H^{1}(\R),
\frac{d}{dt} \scd{x(t)}{\varphi}
+ \scd{x(t)}{i\A' \varphi}=0 , \\
\label{eqweaksolxi}
\forall \varphi\in \D(\A), \ t\mapsto \scd{\xi(t)}{\varphi}\in H^{1}(\R),
\frac{d}{dt} \scd{\xi(t)}{\varphi}
+ \scd{\xi(t)}{i\A \varphi} =  \scd{u(t)}{\CC \varphi} 
.
\end{gather}

The following dual notions of observability and controllability 
are equivalent (cf.~\cite{DR77}).
\begin{defn}
\label{def:exactobs} 
The system (\ref{eqsystx}) 
is {\em exactly observable} in time $T$ at cost $\kappa_{T}$
if the following observation inequality holds:
\begin{gather} 
\label{eqexactobsx}
\forall x_{0}\in X,
\quad
\norm{x_{0}}^{2}\leq \kappa_{T}
\int_{0}^{T}\norm{y(t)}^{2}dt  
. \end{gather}

The system (\ref{eqsystxi}) 
is {\em exactly controllable} in time $T$ at cost $\kappa_{T}$ if 
for all $\xi_{0}$ in $X'$, 
there is a $u$ in $L^{2}(\R;Y')$ such that 
$u(t)=0$ for $t\notin [0,T]$, $\xi(T)=0$ and:
\begin{gather}
\label{eqexactconxi}
\int_{0}^{T}\norm{u(t)}^{2}dt\leq \kappa_{T}\norm{\xi_{0}}^{2}
. \end{gather} 

The {\em controllability cost} for (\ref{eqfsystcon}) in time $T$
is the smallest constant in (\ref{eqexactconxi}),
or in (\ref{eqexactobsx}), still denoted $\kappa_{T}$.
\end{defn}

In this setting, smoother data 
can be controlled by smoother input functions.
The Sobolev space $H^{1}_{0}(0,T;Y')$
is endowed with the homogeneous norm defined by 
$\norm{u}_{1}^{2}=\int_{0}^{T}\norm{\frac{d}{dt}\left(
e^{-i\beta t}u(t)\right)}^{2}dt$,
and its dual is $H^{-1}(0,T;Y)$ with dual norm $\nom{\cdot}$.
Integrating by parts, for all $x_{0}\in X_{-1}$, 
$y(t)=\CC e^{it\A}x_{0}$ satisfies:
\begin{gather*}
\nom{y}=\inf_{\phi \in H^{1}_{0}(0,T;Y')} \left| \int_{0}^{T}
\scd{\CC e^{it(\A-\beta)}(\A-\beta)^{-1}x_{0}}
{\frac{d}{dt}\left(e^{-i\beta t}\phi(t)\right)}dt
\right| / \noo{\phi} .
\end{gather*}
With this remark (and the usual duality argument) we obtain:
\begin{lem}
\label{lem:smoothcontrol}
Let $\beta_{*}=\sup_{t\in [0,T]}\abs{e^{-i\beta t}}^{-2}$
and $\beta^{*}=\sup_{t\in [0,T]}\abs{e^{-i\beta t}}^{2}$.
The admissibility assumptions (\ref{eqadmobs}) and (\ref{eqadmcon}) imply :
$\forall x_{0}\in X$, 
$\nom{\CC e^{it\A}x_{0}}^{2} \leq \beta^{*}K_{T}\nom{x_{0}}^{2}$,
and $\forall u\in H^{1}_{0}(\R;Y')$, 
$\noo{\int_{0}^{T}e^{it\A'}\B u(t)dt}^{2}\leq \beta^{*}K_{T}\noo{u}^{2}$.
Definition~\ref{def:exactobs} implies~:
$\forall x_{0}\in X_{-1}$,
$\nom{x_{0}}^{2}\leq \beta_{*}\kappa_{T}\nom{y}^{2}$,
and equivalently:
for all $\xi_{0}$ in $X'_{1}$, 
there is a $u$ in $H^{1}(\R;Y')$ such that 
$u(t)=0$ for $t\notin (0,T)$, $\xi(T)=0$ and
$\noo{u}^{2}\leq \beta_{*}\kappa_{T}\noo{\xi_{0}}^{2}$.
\end{lem}

The first order control system (\ref{eqfsystcon}) 
and its dual observation system:
\begin{gather}\label{eqfsystobs}
\df(t)-iA f(t)=0, \quad 
f(0)=f_{0}\in H_{1}, \quad
y(t)=Cf(t) 
, \end{gather}
fit into the present setting: 
$X=H_{1}$, $X'=H_{-1}$, $\A$ is $A$ with $\D(\A)=H_{3}$,
$\A'$ is $A'$ with $\D(\A')=H_{1}$,
$\beta=0$, $\beta_{*}=\beta^{*}=1$,
$\CC$ is the $\A$-bounded operator $C$
with $\D(\CC)=H_{2}$ invariant by $(e^{it\A})_{t\geq 0}$.
We shall now explain how the second order control system (\ref{eqsystcon})
and its  dual observation system:
\begin{gather}\label{eqsystobs}
\ddz(t)+Az(t)=0, \quad 
z(0)=z_{0}\in H_{1},\quad \dz(0)=z_{1}\in H_{0},  \quad
y(t)=Cz(t) 
, \end{gather}
also fit into the present setting.

The states $x(t)$ and $\xi(t)$ 
of the systems (\ref{eqsystobs}) and (\ref{eqsystcon}) at time $t$
and their state spaces $X$ and $X'$ are defined by:
\begin{gather*}
\label{eqstate}
x(t)=(z(t),\dz(t))\in X=H_{1}\times H_{0}, 
\quad \xi(t)=(\zeta(t),\dza(t))\in X'=H_{0}\times H_{-1}
. \end{gather*} 
$X$ is a Hilbert space with the ``energy norm'' 
defined by
$\norm{(z_{0},z_{1})}^{2}=\noz{\sA z_{0}}^{2}+\noz{z_{1}}^{2}$,
$X'$ is a Hilbert space with norm 
defined by
$\norm{(\zeta_{0},\zeta_{1})}_{\prime}^{2}
=\noz{\zeta_{0}}^{2}+\nom{\zeta_{1}}^{2}$,
and $X$ is densely continuously embedded in $X'$.
These spaces are dual with respect to the pairing
$\scd{(\zeta_{0},\zeta_{1})}{(z_{0},z_{1})}
=\scz{A^{-1/2}\zeta_{1}}{A^{1/2} z_{0}}-\scz{\zeta_{0}}{z_{1}}$.

The dual second-order systems (\ref{eqsystobs}) and (\ref{eqsystcon}) 
rewrite as dual first order systems (\ref{eqsystx}) and (\ref{eqexactconxi}),
where $u=v$,
$\A $ is defined on the domain 
$\D(\A )=\D(A)\times D(\sA)$ by $\A (z_{0},z_{1})=-i(z_{1},-Az_{0})$,
$\A' $ is an extension of $\A$ to $X'$ with domain $X$,
$\beta=0$, $\beta_{*}=\beta^{*}=1$,
$X_{1}$ is $H_{2}\times H_{1}$ with the norm defined by 
$\norm{(z_{0},z_{1})}^{2}=\norm{\A(z_{0},z_{1})}^{2}
=\noz{\sA z_{1}}^{2}+\noz{Az_{0}}^{2}$,
$\CC\in\L(X_{1},Y)$ is defined by $\CC (z_{0},z_{1})=Cz_{0}$
and $\B \in\L(Y,X'_{-1})$ is defined by $\B y=(0,By)$.

The following admissibility assumptions are then equivalent:
(\ref{eqadm}) for $B$, (\ref{eqadmcon}) for $\B$, (\ref{eqadmobs}) for $\CC$, 
and the admissibility of $C$ for (\ref{eqsystobs}), 
i.e.: $\forall T>0$, $\exists K_{2,T}>0$,
\begin{gather}
\label{eqadmobsz}
\forall x_{0}=(z_{0},z_{1})\in \D(\A), \quad 
\int_{0}^{T}\norm{C z(t)}^{2}dt \leq
K_{2,T}(\noo{z_{0}}^{2}+\noz{z_{1}}^{2})
. 
\end{gather}
In particular, 
$\zeta$ is the unique solution of (\ref{eqsystcon}) 
in $C^{0}(\R;H_{0})\cap C^{1}(\R;H_{-1})$ in the following weak sense:
$\zeta(0)=\zeta_{0}$, $\dza(0)=\zeta_{1}$,
for all $\varphi$ in $\D(A)$,
\begin{gather}
\label{eqweakzeta}
t\mapsto \scz{\zeta(t)}{\varphi}\in H^{2}(\R)
\quad\text{and}\quad  
\frac{d^{2}}{dt^{2}} \scz{\zeta(t)}{\varphi}
+ \scz{\zeta(t)}{A\varphi}=\scz{v(t)}{C\varphi}
\ . \end{gather}
The exact controllability for (\ref{eqsystcon}) in definition~\ref{def:con}
is the usual notion for (\ref{eqsystxi}) in definition~\ref{def:exactobs}.
Similarly, the usual notion of observability 
for (\ref{eqsystx}) in definition~\ref{def:exactobs}
yields the following definition for 
the exact observability in time $T$ at cost $\kappa_{T}$  
of the system (\ref{eqsystobs}):
\begin{gather} 
\label{eqexactobs}
\forall z_{0}\in H_{1}, \forall z_{1}\in H_{0},
\quad
\noz{\sA z_{0}}^{2}+\noz{z_{1}}^{2}
\leq 
\kappa_{T}\int_{0}^{T}\norm{C\dz(t)}^{2}dt 
. \end{gather}


\section{Observability resolvent estimate}
\label{sec:res}

In the general setting for conservative control systems
described in section~\ref{sec:setting}, 
we consider the following observability resolvent estimate:
\begin{gather}
\label{eqRes1}
\exists M>0,  \exists m>0,
\forall x\in\D(\A), \forall \lambda\in \R, 
\quad 
\norm{x}^{2}
\leq M\norm{(\A-\lambda) x}^{2} + m\norm{\CC x}^{2}
. \end{gather}

\begin{thm}
\label{th:selfadj}
The system (\ref{eqsystx}) is exactly observable 
if and only if  the observability resolvent estimate (\ref{eqRes1}) holds.
More precisely, for all $\eps>0$ there is a $C_{\eps}>0$
such that (\ref{eqRes1}) implies (\ref{eqexactobsx})
for all $T> \sqrt{M(\pi^{2}+\eps)}$ 
with $\kappa_{T}=C_{\eps}mT/(T^{2}-M(\pi^{2}+\eps))$.
\end{thm}

We begin by proving two lemmas 
which do not rely on the assumption that $A$ is self-adjoint. 

\begin{lem}
\label{lem:imply}
For all $T>0$, $x_{0}\in \D(\A)$, 
$\lambda\in \R$:
\begin{gather}
\label{eqimply}
\int_{0}^{T} \norm{\CC e^{it\A}x_{0}}^{2}dt
\leq  2T\norm{\CC x_{0}}^{2}
+T^{2}\int_{0}^{T} \norm{\CC e^{it\A}(\A-\lambda)x_{0}}^{2}dt
. \end{gather}
In particular, if the system (\ref{eqsystx}) is exactly observable 
then (\ref{eqRes1}) holds.
\end{lem}
\begin{proof}
Set $x(t)=e^{it\A}x_{0}$, $z(t)=x(t)-e^{it\lambda}x_{0}$ 
and $f=i(\A-\lambda)x_{0}$.
Since $\dx(t)=i\A x(t)=e^{it\A}(i\lambda x_{0}+f)=i\lambda x(t)+e^{it\A}f$,  
we have $\dz(t)=i\lambda z(t)+e^{it\A}f$ 
and therefore $z(t)=\int_{0}^{t}e^{i(t-s)\lambda}e^{is\A}f \, ds$.
We plug it in $x(t)=e^{it\lambda}x_{0}+z(t)$ to estimate:
\begin{gather}
\int_{0}^{T} \norm{\CC x(t)}^{2} dt 
\leq 2\int_{0}^{T} \abs{e^{it\lambda}}^{2}dt \norm{\CC x_{0}}^{2}   
+ 2\int_{0}^{T} t \int_{0}^{t} 
\abs{e^{i(t-s)\lambda}}^{2}\norm{\CC e^{is\A}f}^{2} \, ds\, dt
. \end{gather} 
Since $\lambda\in \R$, 
we have $\abs{e^{it\lambda}}=\abs{e^{i(t-s)\lambda}}^{2}=1$.
Now the inequality: 
\begin{gather*}
\int_{0}^{T} t \int_{0}^{t}F(s) \, ds\, dt
\leq \int_{0}^{T} t \int_{0}^{T}F(s) \, ds\, dt
=(T^{2}/2)\int_{0}^{T}F(s) \, ds
\end{gather*}
with $F(s)=\norm{\CC e^{is\A}f}^{2}$ 
completes the proof of (\ref{eqimply}).

The second statement of lemma~\ref{lem:imply} results from 
applying (\ref{eqadmobs}) and (\ref{eqexactobsx}) 
to (\ref{eqimply}):
it yields (\ref{eqRes1}) with $M=T^{2}\kappa_{T}K_{T}$
and $m=2T\kappa_{T}$. 
\end{proof}

\begin{lem}
\label{lem:conv}
If (\ref{eqRes1}) holds then for all 
$\chi\in C^{1}\comp(\R)$:
\begin{gather}
\label{eqconv}
\forall x_{0}\in X, 
\int \norm{e^{it\A}x_{0}}^{2} \left(\chi^{2}(t)-M\dchi^{2}(t)\right)dt 
\leq 
m\int \norm{\CC e^{it\A}x_{0}}^{2}\chi^{2}(t)dt  
. \end{gather}
\end{lem}
\begin{proof}
Let $x_{0}\in\D(A)$. 
Set $x(t)=e^{it\A}x_{0}$, $z=\chi x$ and $f=\dz -i\A z$.
Since $\dx-i\A x =0$,  we have $f=\dchi x$.
The Fourier transform of $f$ with respect to time is
$\hat{f}(\tau)=(-i\tau-i\A)\hat{z}(\tau)$.
Applying (\ref{eqRes1}) to $\hat{z}(\tau)$, integrating in time,
and the unitarity of the Fourier transform yield:
\begin{gather}
\int \norm{z(t)}^{2} dt 
\leq M\int \norm{f(t)}^{2} dt  + m\int \norm{\CC z(t)}^{2} dt 
. \end{gather} 
Subtracting the first term of the right hand side 
and the density of $\D(A)$ in $X$ complete 
the proof of (\ref{eqconv}). 
\end{proof}

\begin{proof}[Proof of theorem~\ref{th:selfadj}]
The implication is the second part of lemma~\ref{lem:imply}.
The converse results from lemma~\ref{lem:conv}
and the following remark (as in~\cite{BZbb}).

Taking $\chi(t)=\phi(t/T)$ with $\phi\in C^{\infty}\comp(]0,1[)$,
we have
\begin{gather}
\int \norm{\CC e^{it\A}x_{0}}^{2}\chi^{2}(t)dt  
\leq \norm{\phi}_{L^{\infty}}^{2}\int_{0}^{T} \norm{\CC e^{it\A}x_{0}}^{2}dt
\end{gather}
and, since $(e^{it\A})_{t\geq 0}$ is assumed to be a unitary group: 
\begin{gather}
\int \norm{e^{it\A}x_{0}}^{2} \left(\chi^{2}(t)-M\dchi^{2}(t)\right)dt 
=\norm{x_{0}}^{2} I_{T}
\end{gather}
with 
\begin{gather}
I_{T}
=\int\left(\phi^{2}(\frac{t}{T})-\frac{M}{T^{2}}\dphi^{2}(\frac{t}{T})\right)dt
=T \int \phi^{2}(t)dt - \frac{M}{T}\int \dphi^{2}(t)dt
. \end{gather}
For $\phi\neq 0$ and $T$ large enough, 
$I_{T}>0$ so that (\ref{eqconv}) implies (\ref{eqexactobs})
with $\kappa_{T}=m\norm{\phi}_{L^{\infty}}^{2}/I_{T}$. 
In particular, since 
\begin{gather*}
\kappa_{T}=mT\frac{\norm{\phi}_{L^{\infty}}^{2}}{\int \phi^{2}(t)dt}
\left( T^{2}- M\frac{\int \dphi^{2}(t)dt}{\int \phi^{2}(t)dt} 
\right)^{-1}
\text{ and } 
\inf_{\phi\in C^{\infty}\comp(]0,1[)}
\frac{\int \dphi^{2}(t)dt}{\int \phi^{2}(t)dt}=\pi^{2},
\end{gather*} 
for all $\eps>0$, there is a $\phi_{\eps}\in C^{\infty}\comp(]0,1[)$
such that $T> M(\pi^{2}+\eps)$ 
implies $\kappa_{T}=C_{\eps}mT/(T^{2}-M(\pi^{2}+\eps))$ 
with $C_{\eps}=\norm{\phi_{\eps}}_{L^{\infty}}^{2}/\int \phi_{\eps}^{2}(t)dt$.
\end{proof}

\begin{rem}
Observability resolvent estimates like (\ref{eqRes1})
are introduced in~\cite{BZbb}
as sufficient conditions for exact observability.
Theorem~\ref{th:selfadj} for $\CC$ bounded on $X$ is proved in \cite{ZY97},
using a more involved strategy of Liu in~\cite{Liu97}
which our proof shortcuts.
Liu had proved that, 
for a conservative first-order systems with bounded control operator,
exact controllability is equivalent to exponential stability.
From this equivalence 
and the Huang-Pr{\"u}ss condition for exponential stability, 
he deduced an observability resolvent condition 
for conservative second-order systems with bounded observation operator
which he called a frequency domain inequality.
\end{rem}


\section{Fast smoothing controllability}
\label{sec:smooth}

In this section, 
as a substitute to the smoothing effect of holomorphic semigroup 
(used in~\cite{LMctm}),
we introduce the notion of smoothing effect controllability.
More precisely, 
in the general setting for conservative control systems
described in section~\ref{sec:setting}, we prove that 
controllability to a subspace with finite spectrum
and a power-like bound on the cost of fast controls
is implied by the following observability resolvent estimates
(a stronger form of (\ref{eqRes1})):
\begin{gather}
\label{eqresass}
\begin{split}
&\exists m>0, \exists \eps\in]0,1[, 
\exists M:\R\to (0,+\infty), 
\limsup_{\abs{\lambda}\to \infty}\abs{\lambda}^{\eps}M(\lambda)<\infty
\text{ such that:}
\\
& 
\forall x\in \D(\A), \forall \lambda\in \spec(\A), 
\quad 
\norm{x}^{2}
\leq  M(\lambda)\norm{(\A-\lambda) x}^{2} + m\norm{\CC x}^{2}. 
\end{split}
\end{gather}

\begin{thm}
\label{th:ressmoothing}
Assume that $\A$ satisfies (\ref{eqresass}).
$\exists \kappa>0$, $\exists d>0$, 
$\forall T\in ]0,1]$,
$\forall \xi_{0}\in X'$, 
$\exists u\in L^{2}(\mathbb{R};Y)$
such that the solution $\xi\in C^{0}([0,\infty);X')$ of
$$
\dxi(t)-i\A'\xi (t)=Bu(t), \quad \xi(0)=\xi_{0}, 
$$
satisfies $\cd{1}_{\abs{\A'}^{\eps}\geq d/T^{2}}\xi(T)=0$
and 
$\int_{0}^{T}\norm{u(t)}^{2}dt\leq 
\frac{\kappa}{T}\norm{\cd{1}_{\abs{\A'}^{\eps}\geq d/T^{2}}\xi_{0}}^{2}$.
In particular, for all positive $s$, 
$\xi(T)\in \D(\As)$ and 
$\norm{\As \xi(T)}
\leq (1+\sqrt{\kappa K_{T}/T})(d/T^2)^{s/\eps} \norm{\xi_{0}}$.
\end{thm}

\begin{proof} 
Note that (\ref{eqresass}) still holds for all $\lambda\in \R$.
Replacing $\CC$ by $\sqrt{m}\CC$, we assume that $m=1$ 
without loss of generality.

The second statement of the theorem results from applying the first statement
and (\ref{eqadmcon}) to the integral formula expressing $\xi(T)$
in (\ref{eqmildsol}).

The first statement of the theorem is 
the exact controllability in time $T$ of 
the projection of the data 
on the spectral subspace of $X$ 
with spectrum greater than $(d/T^{2})^{1/\eps}$. 
By duality, it is equivalent to the following 
exact observability of data in this spectral subspace:
$\exists \kappa>0$, $\exists d>0$, $\forall T\in ]0,1]$,
\begin{gather}\
\label{eqhighobs}
\forall x_{0}\in  X' \text{ such that } \cd{1}_{\abs{\A}^{\eps}\geq d/T^{2}}x_{0}=x_{0}, 
\quad
\norm{x_{0}}^{2}\leq \frac{\kappa}{T}
\int_{0}^{T}\norm{\CC e^{it\A}x_{0}}^{2}dt 
\ .
\end{gather}
Let $\chi_{T}$ denote a smooth truncation defined by 
$\chi_{T}(t)=\chi(t/T)$ and $\chi\in C\comp(]0,1[)$.
Set $x(t)=e^{it\A}x_{0}$, $z=\chi x$ and $f=i\dz + \A z$.
Since $i\dx+\A x=0$,  we have $f=i\dchiT x$.
The Fourier transform of $f$ with respect to time is
$\hat{f}(\tau)=(\A -\tau)\hat{z}(\tau)$.
With $x=\hat{z}(\tau)$ and $\lambda=\tau$, 
the inequality in (\ref{eqresass}) writes:
$\norm{\hat{z}(\tau)}^{2}
\leq  M(\tau)\norm{\hat{f}(\tau)}^{2} + \norm{\CC \hat{z}(\tau)}^{2}$.
Applying this inequality for $\tau$ greater than a threshold $\mu>0$
and using the unitarity of the Fourier transform yield:
\begin{gather}
\label{eqmu}
\int \norm{z(t)}^{2} dt 
\leq \sup_{\abs{\tau}\geq \mu}\left|\frac{\tau}{\mu}\right|^{\eps}M(\tau)
\int \norm{f(t)}^{2} dt  + \int \norm{\CC z(t)}^{2} dt 
+ \int\limits_{\abs{\tau}< \mu} \norm{\hat{z}(\tau)}^{2} dt
. \end{gather} 
Setting $\mu=2(d/T^{2})^{1/\eps}$ we have 
$\cd{1}_{2\abs{\A}\geq \mu}x_{0}=x_{0}$.
For $\abs{\tau}< 2\mu$ we have 
$\norm{(\A-\tau)^{-1}\cd{1}_{2\abs{\A}\geq \mu}x_{0}}
\leq \norm{x_{0}} (2\mu-\abs{\tau})^{-1}$,
so that using $i^{-1}\partial_{t}e^{it(\A-\tau)}x_{0}=e^{it(\A -\tau)}(\A -\tau)x_{0}$,
and integrating by parts yield  
$\hat{z}(\tau)
=i^{-1}\int \dchiT(t)e^{it(\A -\tau)}(\A -\tau)^{-1}
\cd{1}_{2\abs{\A}\geq \mu}x_{0} dt
$
and 
$\norm{\hat{z}(\tau)}\leq \norm{\dchi}_{L^{1}}(2\mu-\abs{\tau})^{-1}\norm{x_{0}}$.
Therefore
$\int_{\abs{\tau}\leq \mu}\norm{\hat{z}(\tau)}^{2}dt
\leq \frac{2}{\mu}\norm{\dchi}_{L^{1}}^{2}\norm{x_{0}}^{2}$.
Moreover 
$\int \norm{\CC z(t)}^{2}dt \leq
\norm{\chi}_{L^{\infty}}^{2} \int\norm{\CC x(t)}^{2}dt$,
$\int \norm{z(t)}^{2}dt=T \norm{\chi}_{L^{2}}^{2} \norm{x_{0}}^{2}$
and \linebreak 
$\int \norm{f(t)}^{2}dt=T^{-1} \norm{\dchi}_{L^{2}}^{2} \norm{x_{0}}^{2}$.
Hence (\ref{eqmu}) implies:
\begin{gather*}
\norm{x_{0}}^{2}\left(  
\norm{\chi}_{L^{2}}^{2}
- \frac{\norm{\dchi}_{L^{2}}^{2}}{\mu^{\eps}T^{2}}\sup_{\abs{\tau}\geq \mu}\abs{\tau}^{\eps}M(\tau)
-\frac{2\norm{\dchi}_{L^{1}}^{2} }{\mu T}
\right)
\leq \frac{\norm{\chi}_{L^{\infty}}^{2}}{T} \int\norm{\CC x(t)}^{2}dt
\ .
\end{gather*}
Replacing $\mu$ and $x$ by their values, there is a $\kappa'$ 
depending on $\chi$ and $\eps$ such that: 
\begin{gather*}
\norm{x_{0}}^{2}\left(  
1 
- \frac{\kappa'}{d}\sup_{\abs{\tau/2}^{\eps}\geq d/T}\abs{\tau}^{\eps}M(\tau)
-\frac{\kappa' T^{\frac{2}{\eps}-1}}{d^{1/\eps}} 
\right)
\leq \frac{\kappa'}{T} \int\norm{\CC e^{it\A}x_{0}}^{2}dt
\ .
\end{gather*}
Since $\limsup_{\abs{\lambda}\to \infty}\abs{\lambda}^{\eps}M(\lambda)<\infty$,
$\frac{2}{\eps}-1>0$ and $T<1$, taking $d$ large enough independently of $T$
yields a $\kappa>0$ such that (\ref{eqhighobs}) holds.
\end{proof}

\begin{proof}[Proof of theorem~\ref{th:smoothing}]
Since the system (\ref{eqsystcon}) is exactly controllable,
theorem~\ref{th:selfadj} implies 
the corresponding observability resolvent estimate (\ref{eqRes1}), i.e.:
\begin{gather*}
\label{eqRes}
\begin{split}
&
\exists M>0,  \exists m>0,
\forall z_{0}\in H_{2}, \forall z_{1}\in H_{1}, \forall \lambda\in \R, \quad \\
&\noz{\sA z_{0}}^{2}+\noz{z_{1}}^{2}
\leq M\left(
\noz{\sA (-iz_{1}-\lambda z_{0})}^{2}+\noz{iA z_{0}-\lambda z_{1}}^{2}
\right)
+ m\norm{C z_{0}}^{2}
. 
\end{split}
\end{gather*}
For $\lambda\neq 0$ and $z_{1}=i\lambda^{-1}A z_{0}$, 
this estimate writes:
\begin{gather*}
\forall z_{0}\in H_{3}, 
\forall \lambda\in \R^{*}, 
\noz{ \sA z_{0} }^{2}+\frac{1}{\abs{\lambda}^{2}} \noz{Az_{0}}^{2}
\leq \frac{M}{\abs{\lambda}^{2}} \noz{ \sA (A-\lambda^{2}) z_{0} }^{2}
+ m\norm{C z_{0} }^{2}
. \end{gather*}
In particular, since $A$ is positive:  
\begin{gather*}
\forall z_{0}\in H_{3}, \forall \tau\in\spec(A), \quad
\noo{z_{0}}^{2}
\leq \frac{M}{\abs{\tau}} \noo{(A-\tau) z_{0}}^{2}
+ m\norm{C z_{0}}^{2}
. \end{gather*}
Hence the observability resolvent (\ref{eqresass}) 
corresponding to the system (\ref{eqfsystcon})
holds with $M(\lambda)=M/\abs{\lambda}$ and $\eps=1$.
Applying theorem~\ref{th:ressmoothing} with $s=p/2$
completes the proof of theorem~\ref{th:smoothing}.
\end{proof}


\section{Tensor product with a conservative system}
\label{sec:tensor}

Theorem~\ref{th:prod} results from theorem~\ref{th} and the following lemma. 
This trivial lemma is of independent interest.
It says that the controllability cost of a system 
is not changed by taking its tensor product with a conservative system.
It simplifies greatly and improves on
previous results concerning conservative systems distributed in rectangles
(or other product spaces like cylinders or parallelepipeds):
boundary controllability from one whole side (cf. \cite{KLS85})
and semi-internal controllability (cf. \cite{Har89}).
Some applications are given in section~\ref{sec:appli} and~\cite{LMschrocost}.

\begin{lem}
\label{lem:abstprod}
Let $X $, $\Xt$ and $Y$ be Hilbert spaces
and $I$ denote the identity operator on each of them. 
Let $\A :D(\A )\to X $ and $\AAt:D(\AAt)\to \Xt$ be generators of
strongly continuous semigroups of bounded operators on $X $ and $\Xt$.
Let $\CC :D(\CC )\to Y$ be a densely defined operator on $X $
such that $e^{t\A }D(\CC )\subset D(\CC )$ for all $t>0$.
Let $X \ctensprod \Xt$ and $Y\ctensprod \Xt$ 
denote the closure of 
the algebraic tensor products $X \otensprod \Xt$ and $Y\otensprod \Xt$ 
for the natural Hilbert norms.
The operator $\CC \otensprod I:D(\CC )\otensprod \Xt \to Y\ctensprod \Xt$ 
is densely defined on $X \ctensprod \Xt$.

i) The operator $\A \otensprod I+I\otensprod\AAt$ defined on 
the algebraic $D(\A )\otensprod D(\AAt)$ is closable
and its closure, denoted $\A +\AAt$, 
generates a strongly continuous semigroup of bounded operators on $X \ctensprod \Xt$
satisfying:
\begin{gather}
\label{eqnormprod}
\forall t\geq 0, \forall (x,\xt)\in D(\CC )\times \Xt,
\norm{(\CC \otensprod I) e^{t(\A +\AAt)}(x\otensprod \xt)}=\norm{\CC e^{t\A }x}\, \norm{e^{t\AAt}\xt}
\end{gather}

ii) If $i\AAt$ is self-adjoint, then for all $T\geq 0$: 
\begin{gather}
\label{eqcostprod}
\inf_{z\in X \ctensprod \Xt, \norm{z}=1} \int_{0}^{T}\norm{(\CC \otensprod I)e^{t(\A +\AAt)}z}^{2} dt
= \inf_{x\in X , \norm{x}=1} \int_{0}^{T}\norm{\CC e^{t\A }x}^{2} dt
\ .
\end{gather}
\end{lem}
\begin{rem}
When $\CC $ is an admissible observation operator, 
(\ref{eqcostprod}) says that 
the cost of observing $t\mapsto e^{t(\A +\AAt)}$ through $\CC \otensprod I$ in time $T$
is exactly 
the cost of observing $t\mapsto e^{t\A }$ through $\CC $ in time $T$.
The proof of part~i) of lemma~\ref{lem:abstprod} is still valid 
if $X $, $\Xt$ and $Y$ are Banach spaces
and $X \ctensprod \Xt$ and $Y\ctensprod \Xt$ are closures 
with respect to some uniform cross norms (cf. \cite{Sch50}).
\end{rem}

\begin{proof}
Let $G$ denote the generator of 
the strongly continuous semigroup $t\mapsto e^{t\A }\otensprod e^{t\AAt}$
(defined since the natural Hilbert norm is a uniform cross norm, cf. \cite{Sch50}).
Since $D(\A )\otensprod D(\AAt)$ is dense in $X \otensprod \Xt$ and invariant by $t\mapsto e^{tG}$,
it is a core for $G$ (cf. theorem~X.49 in~\cite{RS}).
Since $\A \otensprod I+I\otensprod\AAt=G\res{D(\A )\otensprod D(\AAt)}$, 
it is closable and $\A +\AAt=G$.
Therefore $e^{t(\A+\AAt)}=e^{t\A}\otensprod e^{t\AAt}$ and (\ref{eqnormprod}) follows 
(by the cross norm property).

To prove point~ii), we denote the left and right hand sides of (\ref{eqcostprod})
by $\mathcal{I}_{\A+\AAt}$ and $\mathcal{I}_{\A}$.
Taking $z=x\otensprod\xt$ with $\norm{\xt}=1$, 
$\mathcal{I}_{\A+\AAt}\leq\mathcal{I}_{\A}$ results from (\ref{eqnormprod}).
To prove $\mathcal{I}_{\A+\AAt}\geq\mathcal{I}_{\A}$, 
we only consider the case in which both $X $ and $\Xt$ are 
infinite dimensional and separable
(this simplifies the notation and the other cases are similar).
Let $(e_{n})_{n\in \N}$ and $(\et_{n})_{n\in \N}$ be orthonormal bases for $X $ and $\Xt$.
Since $(e_{n}\otensprod \et_{m})_{n,m\in \N}$ is an orthonormal base for $X \ctensprod \Xt$,
any $z\in X \ctensprod \Xt$ writes: 
\begin{gather*}
z=\sum_{m}x_{m}\otensprod \et_{m} 
\quad \text{with }\ x_{m}=\sum_{n}c_{n,m}e_{n}
\ \text{ and }\ 
\norm{z}^{2}=\sum_{n,m}\abs{c_{n,m}}^{2}=\sum_{m}\norm{x_{m}}^{2}
\ .
\end{gather*}
Since $i\AAt$ is self-adjoint, 
$t\mapsto e^{t\AAt}$ is unitary for all $t\geq 0$
so that  $(e^{t\AAt}\et_{n})_{n\in \N}$ is orthonormal.
Therefore, using (\ref{eqnormprod}):
\begin{gather*}
\norm{\CC e^{t(\A+\AAt)}z}^{2}=\norm{\sum_{m}(\CC e^{t\A}x_{m})\otensprod(e^{t\AAt}\et_{m})}^{2}
=\sum_{m}\norm{\CC e^{t\A}x_{m}}^{2}
\ .
\end{gather*}
By definition, 
$\int_{0}^{T}\norm{\CC e^{t\A}x_{m}}^{2} dt\geq \mathcal{I}_{\A}\norm{x_{m}}^{2}$.
Summing up over $m\in\N$, we obtain:
\begin{gather*}
\int_{0}^{T}\norm{(\CC \otensprod I)e^{t(\A+\AAt)}z}^{2} dt
=\int_{0}^{T}\sum_{m}\norm{\CC e^{t\A}x_{m}}^{2}
\geq \mathcal{I}_{\A}\sum_{m}\norm{x_{m}}^{2}
= \mathcal{I}_{\A}\norm{z}^{2}
\ .
\end{gather*}
This proves $\mathcal{I}_{\A+\AAt}\geq\mathcal{I}_{\A}$ 
and completes the proof of lemma~\ref{lem:abstprod}.
\end{proof}


\section{The fundamental controlled solution}
\label{sec:fcs}

In this section we use theorem~\ref{th:1d} to 
construct a ``fundamental controlled solution'' $k$
of the Schr{\"o}dinger equation on a segment controlled 
by Dirichlet conditions at both ends.

The following proposition 
shows that the upper bound for the controllability cost 
of the Schr{\"o}dinger equation on the segment $[0,L]$ 
controlled at one end
is the same as the controllability cost 
of the Schr{\"o}dinger equation on the twofold segment $[-L,L]$ 
controlled at both ends.

\begin{prop}
\label{prop:twofold}
For any $\alpha > \alpha_{*}$ (cf. definition~\ref{defin:alpha}),
there exists $\gamma >0$ such that, 
for all $L>0$, $T\in\, ]0,\inf(\pi/2,L)^{2}]$
and $\phi_{0}\in H^{-1}(-L,L)$,
there are $g_{-}$ and $g_{+}$ in $L^{2}(0,T)$
such that the solution $\phi\in C^{0}([0,T];H^{-1}(-L,L))$
of the following Schr{\"o}dinger equation on $\left[-L,L\right]$
controlled 
by $g_{-}$ and $g_{+}$:
\begin{equation} 
\label{eqHeattwofold}
\partial_{t}\phi + i\partial_{s}^{2} \phi=0
\quad {\rm in}\ ]0,T[\times ]-L,L[ ,\quad 
\phi\res{s=\pm L} =g_{\pm} ,\quad
\phi\res{t=0} = \phi_{0} 
\end{equation}
satisfies $\phi=0$ at $t=T$ and 
$\displaystyle
\int_{0}^{T}\abs{g_{\pm}(t)}^{2}dt 
\leq \gamma e^{\alpha L^{2}/T }\norm{\phi_{0}}_{H^{-1}(-L,L)}^{2}$.
\end{prop}

\begin{proof}
By duality (cf.~\cite{DR77}), 
it is enough to prove the observation inequality:
$\displaystyle
\exists \gamma >0, 
\forall \phi_{0}\in H^{1}_{0}(-L,L),
\|\phi_{0}\|_{H^{1}}^{2}
\leq \gamma e^{\alpha L^{2}/T}
\|\partial_{s} e^{it\Delta}\phi_{0}{}_{\rceil s=\pm L}
\|_{L^{2}(0,T)^{2}}^{2}
$,
where $\Delta$ denotes $\partial_{s}^{2}$ on $[-L,L]$
with Dirichlet boundary conditions.
Applying theorem~\ref{th:1d} with $N=0$ to the odd part of $\phi_{0}$ 
and with $N=1$ to the even part of $\phi_{0}$
completes the proof of proposition~\ref{prop:twofold}.
\end{proof}

Expressing the solution of (\ref{eqHeattwofold})
with $\phi_{0}=\delta\in H^{-1}(-L,L)$ 
(the Dirac distribution at the origin)
in terms of $g_{\pm}$ by the integral formula 
and applying proposition~\ref{prop:twofold} 
yields the following family of null-controlled solutions 
(depending on $L>0$ and $T>0$ with a good cost estimate)
which we refer to as fundamental controlled solutions.
\begin{cor}
\label{prop:fundcontrol}
For any $\alpha > \alpha_{*}$ (cf. definition~\ref{defin:alpha}),
there exists $\gamma >0$ such that 
$\forall L>0$, $\forall T\in\, ]0,\inf(\pi/2,L)^{2}]$,
$\exists k\in C^{0}([0,T]; H^{-1}(]-L,L[))$
satisfying:
\begin{align}
&\partial_{t}k + i\partial_{s}^{2}k  =  0
\quad {\rm in }\ \mathcal{D}'(]0,T[\times ]-L,L[)\ , 
\label{eqk1} \\
&k\res{t=0}  =  \delta \quad {\rm and }\quad k\res{t=T}  =  0 \ ,
\label{eqk2} \\
&\int_{0}^{T}\norm{k(t,\cdot)}_{H^{-1}(]-L,L[)}^{2}dt
\leq \gamma e^{\alpha L^{2}/T } \ .
\label{eqk3}
\end{align}
\end{cor}


\section{The transmutation of second-order controls into first-order controls}
\label{sec:transmut}

In this section we perform a transmutation 
of a control for the second-order system (\ref{eqsystcon})
into a control for the first-order system (\ref{eqfsystcon})
(cf. (\ref{eqtrans})),
then combine it with theorem~\ref{th:smoothing} into theorem~\ref{th}.
The control transmutation method 
outlined in section~\ref{sec:intro}
proves theorem~\ref{th} only for smoother data, i.e.~:
\begin{prop}
\label{prop:transmut}
If the system (\ref{eqsystcon}) is exactly controllable 
in times greater than $L_{*}$ (cf. definition~\ref{def:con}),
then $\exists \alpha>0$, $\exists \gamma>0$, 
$\forall L>L_{*}$, $\forall T\in \left]0,\inf(1,L)^{2}\right]$,
$\forall \phi_{0}\in H_{1}$, $\exists u \in L^{2}(0,T;Y)$ such that 
the solution $\phi$ of (\ref{eqfsystcon}) 
satisfies $\phi(T)=0$ and 
$
\int_{0}^{T} \norm{u(t)}^{2}dt
\leq \kappa_{2,L}\gamma e^{\alpha L^{2}/T }  \noo{\phi_{0}}^{2}$,
where $\kappa_{2,L}$ is defined in (\ref{eqkappa2}).
\end{prop}

\begin{proof}
Let $L>L_{*}$. 
Since (\ref{eqsystcon}) is exactly controllable in time $L$,
by lemma~\ref{lem:smoothcontrol} 
(applied to the reduction of (\ref{eqsystcon}) to the first-order setting
described after the statement of this lemma):
for all $\zeta_{0}\in H_{1}$ and $\zeta_{1}\in H_{0}$,
there is a $v\in H^{1}(\R; Y)$ such that 
$v(s)=0$ for $s\notin (0,L)$, 
the solution $\zeta$ of (\ref{eqsystcon})
satisfies $\zeta(L)=\dza(L)=0$ and 
\begin{gather}
\label{eqcostv}
\int \norm{\dv(t)}^{2}dt\leq \kappa_{2,L}
\left( \noo{\zeta_{0}}^{2}+\noz{\zeta_{1}}^{2} \right)
. \end{gather} 

Let $\alpha>\alpha_{*}$ and $T\in]0,\inf(1,L^{2})[$.  
Let $\gamma>0$ and 
$k\in C^{0}([0,T]; H^{-1}(]-L,L[))$ 
be the corresponding constant and fundamental controlled solution
given by corollary~\ref{prop:fundcontrol}.
We define $\kk\in C^{0}([0,\infty); H^{-1}(\R))$ 
as the extension of $k$ by zero,
i.e. $\kk(t,s)=\bar{k}(t,s)$ on $[0,T]\times]-L,L[$
and $\kk$ is zero everywhere else.
It inherits from $k$ the following properties 
\begin{align}
&\partial_{t}\kk + i\partial_{s}^{2} \kk=0
\quad {\rm in }\  \mathcal{D}'(]0,T[\times ]-L,L[)\ , 
\label{eqkk1} \\
&\kk\res{t=0}  =  \delta \quad {\rm and }\quad \kk\res{t=T}  =  0 \ ,
\label{eqkk2} \\
&\int_{0}^{T} \norm{\kk(t,\cdot)}_{H^{-1}(\R)}^{2} dt
\leq \gamma e^{\alpha L^{2}/T } \ .
\label{eqkk3}
\end{align}

Let $\phi_{0}\in H_{1}$ be an initial data for (\ref{eqfsystcon}).
Let $\zeta$ and $v$
be the corresponding solution and control function for (\ref{eqsystcon}) 
with data $\zeta_{0}=\phi_{0}$ and $\zeta_{1}=0$.
We define $\zz\in C^{0}(\R;H_{1})\cap C^{1}(\R;H_{0})$
and $\vv\in H^{1}(\mathbb{R};Y)$ 
as the extensions of $\zeta$ and $v$ 
by reflection with respect to $s=0$,
i.e. $\zz(s)=\zeta(s)=\zz(-s)$ and $\vv(s)=v(s)=\vv(-s)$ for $s\geq 0$.
Since $\zeta_{1}=\zeta(L)=\dza(L)=0$, 
$\zz$ is the unique solution in $C^{0}(\R;H_{0})\cap C^{1}(\R;H_{-1})$ of: 
\begin{gather*}
\ddzz(t)+\Am\zz (t)=B\vv(t), \quad 
\zz(0)=\phi_{0},\ \dzz(0)=0 
, \end{gather*}
in particular in the following weak sense (as in (\ref{eqweakzeta})):
for all $\varphi$ in $H_{2}$, 
\begin{gather} 
\label{eqzz}
s\mapsto \scz{\zz(s)}{\varphi}\in H^{2}(\R)
\text{ and }
\frac{d^{2}}{ds^{2}} \scz{\zz(s)}{\varphi}
+ \scz{\zz(s)}{A\varphi}=\scz{\vv(s)}{C\varphi}
. \end{gather}
Equation (\ref{eqcostv}) implies the following cost estimate for $\vv$:
\begin{gather} 
\label{eqvv}
\int \norm{\dvv(s)}^{2}ds
\leq 2\int \norm{v(s)}^{2}ds
\leq 2\kappa_{2,L}\noo{\phi_{0}}^{2}
. \end{gather}

The main idea of our proof is to use $\kk$
as a kernel to transmute $\zz$ and $\vv$
into a solution $\phi$ and a control $u$ 
for (\ref{eqfsystcon}). 
The transmutation formulas:
\begin{equation} 
\label{eqtrans}
\phi(t)=\int \kk(t,s)\zz(s)\, ds  
 \quad {\rm and } \quad \forall t>0, \
u(t)=-i\int \kk(t,s)\vv(s)\, ds 
\ ,
\end{equation}
define $\phi\in C^{0}([0,\infty);H_{0})$ and $u\in L^{2}([0,\infty);Y)$
since $\kk\in C^{0}([0,\infty); H^{-1}(\R))$, 
$\zz\in H^{1}(\mathbb{R};H_{0})$ and $\vv\in H^{1}(\mathbb{R};Y)$.
The property $(\ref{eqkk2})$ of $\kk$ implies 
$\phi(0)=\phi_{0}$ and $\phi(T)=0$.
Since $\zz(s)=\dot{\zz}(s)=0$ for $|s|=L$,
equations $(\ref{eqzz})$ and $(\ref{eqkk1})$ imply,
by integrating by parts, 
for all $\varphi$ in $H_{3}$:
\begin{equation}
  \label{equ1}
t\mapsto \scz{\phi(t)}{\varphi}\in H^{1}(0,\infty),\ 
\frac{d}{dt} \scz{\phi(t)}{\varphi} + \scz{\phi(t)}{iA\varphi} 
=\scz{u(t)}{C\varphi}
. \end{equation}
This is the equation (\ref{eqweaksolxi}) corresponding to (\ref{eqfsystcon}),
i.e. with the settings described after (\ref{eqfsystobs}).
Therefore $\phi$ and $u$ satisfy  (\ref{eqfintcon}).

Since 
$\int_{0}^{T} \norm{u(t)}^{2}dt\leq 
\int_{0}^{T} \norm{\kk(t,\cdot)}_{H^{-1}(\R)}^{2}dt \int \norm{\dvv(s)}^{2}ds$,
equations (\ref{eqkk3}) and (\ref{eqvv}) 
imply the cost estimate 
which completes the proof of proposition~\ref{prop:transmut}.
\end{proof}

\begin{proof}[Proof of theorem~\ref{th}]
Let $\alpha > \alpha_{*}$, $L>L_{*}$ and $\eps\in ]0,1[$.

According to theorem~\ref{th:smoothing} with $p=2$: 
$\exists \kappa>0$, $\exists d>0$, 
$\forall T\in ]0,1]$, $\forall \phi_{0} \in H_{-1}$, 
$\exists u_{1} \in L^{2}([0,\eps T];Y)$ such that 
the solution $\phi\in C^{0}([0,\eps T];H_{-1})$ of (\ref{eqfsystcon})
with $u=u_{1}$ on $[0,\eps T]$
satisfies $\phi(T)\in H_{1}$,
$\noo{\phi(T)}
\leq \nom{\phi_{0}}(1+\sqrt{\kappa K_{1,\eps T}/\eps T})d/(\eps T)^2$,
and 
$\int_{0}^{T}\norm{u_{1}(t)}^{2}dt\leq\frac{\kappa}{\eps T}\nom{\phi_{0}}^{2}$.
Therefore, according to proposition~\ref{prop:transmut}, 
$\exists \alpha>0$, $\exists \gamma>0$, 
$\forall T\in \left]0,\inf(1,L)^{2}\right]$,
$\forall \phi_{0} \in H_{-1}$, 
$\exists u_{2} \in L^{2}([\eps T,T];Y)$ such that 
the solution $\phi\in C^{0}([0,T];H_{-1})$ of (\ref{eqfsystcon})
with $u=u_{1}$ on $[0,\eps T]$ and $u=u_{2}$ on $[\eps T,T]$
satisfies $\phi(T)=0$ and 
$
\int_{\eps T}^{T} \norm{u_{2}(t)}^{2}dt \leq 
\kappa_{2,L}\gamma e^{\alpha L^{2}/(T-\eps T) }\nod{\phi_{0}}^{2}$.
Since $\int_{0}^{T} \norm{u(t)}^{2}dt=
\int_{0}^{\eps T}\norm{u_{1}(t)}^{2}dt+\int_{\eps T}^{T}\norm{u_{2}(t)}^{2}dt$,
the controllability cost $\kappa_{1,T}$ in definition~\ref{def:con}
satisfies for all $T\in \left]0,\inf(1,L)^{2}\right]$:
\begin{gather*}
\kappa_{1,T}\leq
\frac{\kappa}{\eps T}+
\left(1+\sqrt{\kappa K_{1,1}/\eps T}\right)^{2}\frac{d^{2}}{(\eps T)^4}\kappa_{2,L}\gamma 
\exp\frac{\alpha L^{2}}{(1-\eps)T} 
.
\end{gather*}
Therefore 
$\displaystyle
\limsup_{T\to 0} T \ln \kappa_{1,T} \leq \alpha L^{2}/(1-\eps)$.
Letting $\alpha$, $L$ and $\eps$  tend respectively to 
$\alpha_{*}$, $L_{*}$ and $0$ 
completes the proof of (\ref{equb}).
\end{proof}


\section{Geometric bounds on the cost of fast boundary controls 
for Schr{\"o}dinger equations}
\label{sec:appli}

When the second-order equation (\ref{eqsystcon}) 
has a finite propagation speed
and is controllable, 
the control transmutation method yields 
geometric upper bounds on the cost of fast controls 
for the first-order equation (\ref{eqfsystcon}).
This was illustrated in \cite{LMschrocost}
on the internal controllability of Schr{\"o}dinger equations
on Riemannian manifolds
which have the wave equation as corresponding second-order equation.
Similar lower bounds proved in \cite{LMschrocost}
(without assuming the controllability of the wave equation)
imply that 
the upper bounds are optimal with respect to time dependence.
In this section, we illustrate the control transmutation method 
on the analogous boundary control problem for Schr{\"o}dinger equations.

Let $(M,g)$ be a smooth connected compact 
$n$-dimensional Riemannian manifold with metric $g$ 
and smooth boundary $\partial M$. 
When $\partial M\neq\emptyset$, $M$ denotes the interior 
and $\adh{M}=M\cup\partial{M}$.
Let $\Delta$ denote the (negative) Laplacian on $(M,g)$
and $\partial_{\nu}$ denote the exterior Neumann vector field on $\partial M$.
The characteristic function of a set $S$ is denoted by $\chi_{S}$.

Let $H_{0}=L^{2}(M)$.
Let $A$ be defined by $Af=-\Delta f$ 
on $\D(A)=H^{2}(M)\cap H^{1}_{0}(M)$.
Let $C$ be defined from $\D(A)$ to $Y=L^{2}(\partial M)$
by $Cf=\partial_{\nu} f\res{\Gamma}$
where $\Gamma$ is an open subset of $\partial M$.
With this setting, 
(\ref{eqfsystcon}) is 
a Schr{\"o}dinger equation,
(\ref{eqsystcon}) is a scalar wave equation,
and these equations are controlled by 
the Dirichlet boundary condition on $\Gamma$.
In particular (\ref{eqsystcon}) writes:
\begin{gather}\label{eqwavecon}
\begin{split}
&\partial_{t}^{2} \zeta-\Delta\zeta =0\ \text{ on }\ \R_{t}\times M, 
\quad 
\zeta=\chi_{\Gamma} v\ \text{ on }\ \R_{t}\times \partial M, \\
&\zeta(0)=\zeta_{0}\in L^{2}(M),\ \dza(0)=\zeta_{1}\in H^{-1}(M) ,\ 
v\in L^{2}\loc(\R;L^{2}(\partial M)) 
, \end{split}
\end{gather}
It is well known that $C$ is an admissible observation operator 
for the wave equation (\ref{eqsystobs}) 
and the Schr{\"o}dinger equation (\ref{eqfsystobs}) 
(cf. e.g. corollary~3.9 in \cite{BLR92} and \cite{Leb92}).
To ensure the exact controllability of the wave equation
we use the geometric optics condition of Bardos-Lebeau-Rauch
(specifically example~1 after corollary~4.10 in \cite{BLR92}):
\begin{gather}\label{eqGC}
\text{\parbox{\textequation}{
There is a positive constant  $L_{\Gamma}$ 
such that every generalized geodesic of length greater than $L_{\Gamma}$ 
passes through $\Gamma$ at a non-diffractive point.
}}
\end{gather} 
Generalized geodesics are the rays of geometrical optics 
(we refer to \cite{Mil02} for a presentation of this condition 
with a discussion of its significance).
We make the additional assumption
that they can be uniquely continued at the boundary $\partial M$.
As in~\cite{BLR92}, to ensure this, we may assume either that 
$\partial M$ has no contacts of infinite order with its tangents
(e.g. $\partial M=\emptyset$),
or that $g$ and $\partial M$ are real analytic.
For instance, we recall that (\ref{eqGC}) holds when 
$\Gamma$ contains a closed hemisphere 
of a Euclidean ball $M$ of diameter $L_{\Gamma}/2$,
or when $\Gamma=\partial M$ and $M$ is a strictly convex bounded Euclidean set
which does not contain any segment of length $L_{\Gamma}$. 
\begin{thm}[\cite{BLR92}]
\label{th:BLR}
If (\ref{eqGC}) holds then the wave equation (\ref{eqwavecon}) 
is exactly controllable in any time greater than $L_{\Gamma}$.
\end{thm}

Thanks to this theorem, theorem~\ref{th:prod} implies:
\begin{thm}
\label{th:geomprod}
Let $\Mt$ be a smooth complete $\nt$-dimensional Riemannian manifold
and $\Deltat$ denote the Laplacian on $\Mt$ 
with the Dirichlet boundary condition.
Let $\gamma$ denote the subset $\Gamma\times \Mt$ 
of $\partial(M\times\Mt)$.
If (\ref{eqGC}) holds then the Schr{\"o}dinger equation: 
\begin{gather*}
\begin{split}
&i\partial_{t}\phi-(\Delta+\Deltat)\phi =0\ 
\text{ on }\ \R_{t}\times M\times\Mt, 
\quad 
\phi=\chi_{\gamma} u\ \text{ on }\ \R_{t}\times \partial (M\times\Mt), \\
&\phi(0)=\phi_{0}\in L^{2}(\Mt;H^{-1}(M)),\ 
u\in L^{2}\loc(\R;L^{2}(\partial (M\times\Mt))) 
, \end{split}
\end{gather*}
is exactly controllable in any time $T$ at a cost $\kapt_{T}$
which satisfies the following upper bound 
(with $\alpha_{*}$ as in theorem~\ref{th:1d}):
$
\displaystyle \limsup_{ T\to 0}  
T \ln \kapt_{T} \leq \alpha_{*} L_{\Gamma}^{2}
$.
\end{thm}
\begin{rem}
\label{rem:geomprod}
For $\Mt= \emptyset$, the controllability was proved in \cite{Leb92}.
As in  \cite{Leb92}, this results extends to the plate equation.
The boundary controllability of a rectangular plate from one side 
was proved in \cite{KLS85} (theorem~2).
When $M$ is a segment and $\Mt$ is a line, 
theorem~\ref{th:geomprod} extends this result to an infinite strip.
\end{rem}
\begin{rem}
In particular, theorem~\ref{th:geomprod} shows that 
the geometric optics condition 
is not necessary for the controllability cost of the Schr{\"o}dinger equation
to grow at most like $\exp(C/T)$ as $T$ tends to $0$.
Indeed, any geodesic of $\Mt$ yields a geodesic of $M\times\Mt$
in a slab $\set{x}\times\Mt$ with $x\in M$,
and this geodesic does not pass through the control region $\gamma$
since the slab does not intersect the boundary set $\partial M\times \Mt$.
\end{rem}

\section*{Acknowledgments}
The investigation of observability resolvent estimates
in sections~\ref{sec:res} and~\ref{sec:smooth}
originates in~\cite{BZbb}.
I am thankful to N.~Burq and M.~Zworski for discussions 
on their manuscript.



\providecommand{\bysame}{\leavevmode\hbox to3em{\hrulefill}\thinspace}
\providecommand{\MR}{\relax\ifhmode\unskip\space\fi MR }
\providecommand{\MRhref}[2]{%
  \href{http://www.ams.org/mathscinet-getitem?mr=#1}{#2}
}
\providecommand{\href}[2]{#2}

\end{document}